\documentclass[journal]{IEEEtran}

\input{xhc.sty}

\begin{document}

\title{Optimal Tap Setting of Voltage Regulation Transformers Using Batch Reinforcement Learning}

\author{Hanchen~Xu,~\IEEEmembership{Student Member,~IEEE,}
        Alejandro~D.~Dom\'{i}nguez-Garc\'{i}a,~\IEEEmembership{Member,~IEEE}
        and Peter~W.~Sauer,~\IEEEmembership{Life~Fellow,~IEEE}
\thanks{The authors are with the Department of Electrical and Computer Engineering at the University of Illinois at Urbana-Champaign, Urbana, IL 61801, USA. Email: \{hxu45, aledan, psauer\}@illinois.edu.}
}

\maketitle

\begin{abstract}
In this paper, we address the problem of setting the tap positions of load tap changers (LTCs) for voltage regulation in radial power distribution systems under uncertain load dynamics.
The objective is to find a policy to determine the tap positions that only uses measurements of voltage magnitudes and topology information so as to minimize the voltage deviation across the system.
We formulate this problem as a Markov decision process (MDP), and propose a batch reinforcement learning (RL) algorithm to solve it.
By taking advantage of a linearized power flow model, we propose an effective algorithm to estimate the voltage  magnitudes under different tap settings, which allows the RL algorithm to explore the state and action spaces freely offline without impacting the system operation.
To circumvent the ``curse of dimensionality" resulted from the large state and action spaces, we propose a sequential learning algorithm to learn an action-value function for each LTC, based on which the optimal tap positions can be directly determined.
The effectiveness of the proposed algorithm is validated via numerical simulations on the IEEE 13-bus and 123-bus distribution test feeders.
\end{abstract}



\section{Introduction} \label{sec:intro}

\IEEEPARstart{V}{oltage} regulation transformers---also referred to as load tap changers (LTCs)---are widely utilized in power distribution systems to regulate the voltage magnitudes along a feeder.
Conventionally, the tap position of each LTC is controlled through an automatic voltage regulator based on local voltage measurements \cite{kundur1994power}.
This approach, albeit simple and effective, is not optimal in any sense, and may result in frequent actions of the LTCs, thus, accelerating wear and tear \cite{robbins2016tap}.
Particularly, the voltage deviation may not be minimized.
In the context of transmission systems, transformer tap positions are optimized jointly with active and reactive power generation by solving an optimal power flow (OPF) problem, which is typically cast as a mixed-integer programming problem (see, e.g., \cite{liu1992shunt, salem1997tap} and references therein). 
Similar OPF-based approaches are also adopted in power distribution systems.
For example, in \cite{robbins2016tap}, the authors cast the optimal tap setting problem as a rank-constrained semidefinite program that is further relaxed by dropping the rank-one constraint, which avoids the non-convexity and integer variables, and thus, the problem can be solved efficiently.
OPF-based approaches have also been utilized to determine the optimal reactive power injection from distributed energy resources so as to regulate voltage in a distribution network \cite{zhu2016fast, robbins2016reactive}. 

While these OPF-based approaches are effective in regulating voltages, they require complete system knowledge, including active and reactive power injections, and transmission/distribution line parameters.
While it may be reasonable to assume that such information in available for transmission systems, the situation in distribution systems is quite different.
Accurate line parameters may not be known and power injections at each bus may not be available in real time, which prevents the application of OPF-based approaches \cite{xu2018voltage}. 
In addition, OPF-based approaches typically deal with one snapshot of system conditions, and assume loads remain constant between two consecutive snapshots.
Therefore, the optimal tap setting problem needs to be solved for each snapshot in real time.

In this paper, we develop an algorithm that can find a policy for determining the optimal tap  positions of the LTCs in a power distribution system under uncertain load dynamics without any information on power injections or line parameters; the algorithm requires only voltage magnitude measurements and system topology information.
Specifically, the optimal tap setting problem is cast as a Markov decision process (MDP), which can be solved using reinforcement learning (RL) algorithms.
Yet, adequate state and action samples that sufficiently explore the MDP state and action spaces are needed.
However, it is hard to obtain such samples in real power systems since this requires changing tap settings and other controls to  excite the system and record voltage responses, which may jeopardize system operational reliability and incur economic costs.
To circumvent this issue, we take advantage of a linearized power flow model and develop an effective algorithm to estimate voltage magnitudes under different tap settings so that the state and action spaces can be explored freely offline without impacting the real system.

The dimension of the state and action spaces increases exponentially as the number of LTCs grows, which causes the issue known as the ``curse of dimensionality" and makes the computation of the optimal policy intractable \cite{sutton2018reinforcement}.
To circumvent the ``curse of dimensionality," we propose an efficient batch RL algorithm---the least squares policy iteration (LSPI) based sequential learning algorithm---to learn an action-value function sequentially for each LTC. 
Once the learning of the action-value function is completed, we can determine the policy for optimally setting the LTC taps.
We emphasize that the optimal policy can be computed offline, where most computational burden takes place.
However, when executed online, the required computation to find the optimal tap positions is minimal.
The effectiveness of the proposed algorithm is validated through simulations on two IEEE distribution test feeders.

The remainder of the paper is organized as follows.
Section \ref{sec:prelim} introduces a linearized power flow model that includes the effect of LTCs and describes the optimal tap setting problem.
Section \ref{sec:mdp} provides a primer on MDPs and the LSPI algorithm.
Section \ref{sec:formulation} develops an MDP-based formulation for the optimal tap setting problem and Section \ref{sec:algo} proposes an algorithm to solve this problem.
Numerical simulation results on two IEEE test feeders are presented in Section \ref{sec:simu}.
Concluding remarks are provided in Section \ref{sec:con}.

\section{Preliminaries} \label{sec:prelim}

In this section, we review a linearized power flow model for power distribution systems, and modify it to include the effect of LTCs.
We also describe the LTC tap setting problem.

\subsection{Power Distribution System Model}

Consider a power distribution system that consists of a set of buses indexed by the elements in $\calN = \{0, 1, \cdots, N\}$, and a set of transmission lines indexed by the elements in $\calL = \{1, \cdots, L\}$.
Each line $\ell \in \calL$ is associated with an ordered pair $(i,j) \in \calN \times \calN$.
Assume bus $0$ is an ideal voltage source that corresponds to a substation bus, which is the only connection of the distribution system to the bulk power grid. 

Let $V_i$ denote the magnitude of the voltage at bus $i$, $i \in \calN$, and define $v_i := V_i^2$; note that $u_0$ is a constant since bus $0$ is assumed to be an ideal voltage source.
Let $p_i$ and $q_i$ denote the active power injection and reactive power injection at bus $i$, $i \in \calN$, respectively. 
For each line $\ell \in \calL$ that is associated with $(i,j)$, let $p_{ij}$ and $q_{ij}$ respectively denote active and reactive power flows on line $(i,j)$, which are positive if the flow of power is from bus $i$ to bus $j$ and negative otherwise.
Let $r_\ell$ and $x_\ell$ denote the resistance and reactance of line $\ell$, $\ell \in \calL$.
For a radial power distribution system, the relation between squared voltage magnitudes, power injections, and line power flows, can be captured by the so-called LinDisfFlow model \cite{33bus} as follows:
\begin{subequations} \label{eq:LDF}
\begin{align}
    p_{ij} &= -p_j + \sum_{k:(j,k)\in \calL} p_{jk}, \\
    q_{ij} &= -q_j + \sum_{k:(j,k)\in \calL} q_{jk}, \\
    v_i - v_j &= 2 (r_{\ell} p_{ij} + x_{\ell} q_{ij}), \label{eq:LDF-V}
\end{align}
\end{subequations}
where $\ell$ is associated with $(i, j)$.

Define $\bm{r}=[r_1, \cdots, r_L]^\top$ and $\bm{x}=[x_1, \cdots, x_L]^\top$.
Let $\tdbdM =[\tilde{M}_{i \ell}] \in \real^{(N+1) \times L}$, with $\tilde{M}_{i \ell} = 1$ and $\tilde{M}_{j \ell} = -1$ if line $\ell$ is associated with $(i, j)$, and all other entries equal to zero.
Let $\bm{m}^\top$ denote the first row of $\tdbdM$ and $\bm{M}$ the matrix that results by removing $\bm{m}^\top$ from $\tdbdM$.
For a radial distribution system, $L=N$, and $\bm{M}$ is invertible.
Define $\bm{v} = [v_1, \cdots, v_N]^\top$, $\bm{p} = [p_1, \cdots, p_N]^\top$, and $\bm{q} = [q_1, \cdots, q_N]^\top$.
Then, the LinDistFlow model in \eqref{eq:LDF} can be written as follows:
\begin{align} \label{eq:LDF-Vec}
    \bm{M}^\top \bm{v} + \bm{m} v_0 = 2\diag{\bm{r}} \bm{M}^{-1} \bm{p} + 2\diag{\bm{x}} \bm{M}^{-1} \bm{q}, 
\end{align}
where $\diag{\cdot}$ returns a diagonal matrix with the entries of the argument as its diagonal elements.

\begin{figure}[!t]
\centering
\includegraphics[width=2.2in]{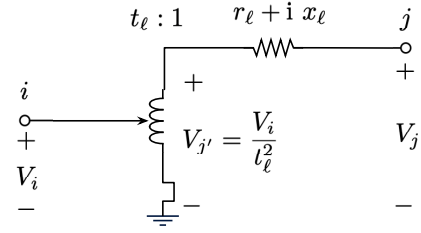}
\caption{Load tap changer model.}
\label{fig:XFM}
\end{figure}

The standard model for an LTC in the literature is shown in Fig. \ref{fig:XFM} (see, e.g., \cite{kundur1994power}), where $\imagi = \sqrt{-1}$, line $\ell$ is associated with $(i,j)$, and $t_\ell$ is the tap ratio of the LTC on line $\ell$.
Typically, the tap ratio can possibly take on $33$ discrete values ranging from $0.9$ to $1.1$, by an increment of $5/8$\% p.u., i.e., $t_\ell \in \mathcal{T} = \{0.9, 0.90625, \cdots, 1.09375, 1.1\}$ \cite{kundur1994power}.
Let $\Delta t_\ell \in \Delta \mathcal{T} = \{0, \pm  0.00625, \cdots, \pm 0.19375, \pm 0.2\}$ denote the set of all feasible LTC tap ratio changes.
We index the $33$ tap positions by $-16, \cdots, -1, 0, 1, \cdots, 16$ for convenience.

Let $\calL^t$ denote the set of lines with LTCs and let $|\calL^t| = L^t$, where $|\cdot|$ denotes the cardinality of a set.
For line $\ell$ that is associated with $(i, j)$, if $\ell \in \calL^t$, the voltage relation in the LinDistFlow model, i.e., \eqref{eq:LDF-V}, needs to be modified as follows:
\begin{equation} \label{eq:LDF-XFM}
    \frac{1}{t_\ell^2} v_i - v_j = 2 (r_\ell p_{ij} + x_\ell q_{ij}).
\end{equation}
Define $\bm{t} = [t_\ell]^\top$ and $\Delta \bm{t} = [\Delta t_\ell]^\top$, $\ell \in \calL^t$.
Let $\tdbdM(\bm{t}) =[\tilde{M}_{i \ell}(\bm{t})] \in \bbR^{(N+1) \times L}$, with $\tilde{M}_{i \ell}(\bm{t}) = 1$ and $\tilde{M}_{j \ell}(\bm{t}) = -1$ if line $\ell \in \calL \setminus \calL^t$, $\tilde{M}_{i \ell}(\bm{t}) = \frac{1}{t_\ell^2}$ and $\tilde{M}_{j \ell}(\bm{t}) = -1$ if line $\ell \in \calL^t$, and all other entries equal to zero.
Let $\bm{m}(\bm{t})^\top$ denote the first row of $\tdbdM(\bm{t})$ and $\bm{M}(\bm{t})$ the matrix that results by removing $\bm{m}(\bm{t})^\top$ from $\tdbdM(\bm{t})$.
The matrix $\bm{M}(\bm{t})$ is non-singular when the power distribution system is connected.
Then, the modified matrix-form LinDistFlow model that takes into account the LTCs is given by:
\begin{align} \label{eq:LDF-Vec-LTC}
    \bm{M}(\bm{t})^\top \bm{v} + \bm{m}(\bm{t}) v_0 = 2\diag{\bm{r}} \bm{M}^{-1} \bm{p} + 2\diag{\bm{x}} \bm{M}^{-1} \bm{q}.
\end{align}

\subsection{Optimal Tap Setting Problem}

To effectively regulate the voltages in a power distribution system, the tap positions of LTCs need to be set appropriately.
The objective of the optimal tap setting problem is to find a policy $\bm{\pi}$ that determines the LTC tap ratio so as to minimize the voltage deviation from some reference value, denoted by $\bm{v}^{\star}$, based on current tap ratios and measurements of the voltage magnitudes, i.e., $\bm{\pi}: (\bm{t},\bm{v}) \rightarrow \Delta \bm{t},~\bm{t} \in \calT^{L^t},~\bm{v} \in \bbR^N, \Delta \bm{t} \in  \Delta \calT^{L^t}$.
Throughout this paper, we make the following two assumptions:
\begin{itemize}
	\item[{\textbf A1}.] The distribution system topology is known but the line parameters are unknown.
	\item[{\textbf A2}.] The active and reactive power injections are not measured and their probability distributions are unknown.
\end{itemize}

\section{Markov Decision Process and Batch Reinforcement Learning} \label{sec:mdp}

In this section, we provide some background on MDPs and the batch RL algorithm, a type of data efficient and stable algorithm for solving MDPs with unknown models.

\subsection{Markov Decision Process}

An MDP is defined as a 5-tuple $(\calS, \calA, \calP, \calR, \gamma)$, where $\calS$ is a finite set of states, $\calA$ is a finite set of actions, $\calP$ is a Markovian transition model that denotes the probability of transitioning from one state into another after taking an action, $\calR: \calS \times \calA \times \calS \rightarrow \bbR$ is a reward function such that, for $\bm{s}, \bm{s}' \in \calS$ and $\bm{a} \in \calA$, $r = \calR(\bm{s}, \bm{a}, \bm{s}')$ is the reward obtained when the system transitions from state $\bm{s}$ into state $\bm{s}'$ after taking action $\bm{a}$, and $\gamma \in [0, 1)$ is a discount factor (see, e.g., \cite{lagoudakis2003least}).\footnote{These definitions can be directly extended to the case where the the set of states is infinite. Due to space limitation, this case is not discussed in detail here.}
We refer to the 4-tuple $(\bm{s}, \bm{a}, r, \bm{s}')$, where $\bm{s}'$ is the state following $\bm{s}$ after taking action $\bm{a}$ and $r=\calR(s, a, s')$, as a transition.

Let $\bm{S}_k$ and $\bm{A}_k$ denote the state and action at time instant $k$, respectively, and $R_k$ the reward received after taking action $\bm{A}_k$ in state $\bm{S}_k$.
Let $\bbP$ denote the probability operator; then, $\calP_k(\bm{s}' | \bm{s}, \bm{a}) := \prob{\bm{S}_{k+1} = \bm{s}' | \bm{S}_k = \bm{s}, \bm{A}_k = \bm{a}}$ is the probability of transitioning from state $\bm{s}$ into state $\bm{s}'$ after taking action $\bm{a}$ at instant $k$.
Throughout this paper, we assume time-homogeneous transition probabilities, hence we drop the subindex $k$ and just write $\calP(s'|s, a)$.

Let $\bar{R}: \calS \times \calA \rightarrow \bbR$ denote the expected reward for a state-action pair $(\bm{s}, \bm{a})$; then, we have
\begin{align} \label{eq:expected_reward}
	\bar{R}(\bm{s}, \bm{a}) = \expect{R} = \sum \limits_{\bm{s'} \in \calS} \calR(\bm{s}, \bm{a}, \bm{s}') \calP(\bm{s}' | \bm{s}, \bm{a}),
\end{align}
where $\expect{\cdot}$ denotes the expectation operation.
The total discounted reward from time instant $k$ and onwards, denoted by $G_k$, also referred to as the return, is given by
\begin{align} \label{eq:return}
	G_k = \sum_{k'=k}^\infty \gamma^{k'-k} R_{k'}.
\end{align}

A deterministic policy $\bm{\pi}$ is a mapping from $\calS$ to $\calA$, i.e., $\bm{a} = \bm{\pi}(\bm{s}), \bm{s} \in \calS, \bm{a} \in \calA$.
The action-value function under policy $\bm{\pi}$ is defined as follows:
\begin{align} \label{eq:Q_def}
	Q^{\bm{\pi}}(\bm{s}, \bm{a}) = \expect{G_k | \bm{S}_k = \bm{s}, \bm{A}_k = \bm{a}; \bm{\pi}},
\end{align}
which is the expected return when taking action $\bm{a}$ in state $\bm{s}$, and following policy $\bm{\pi}$ afterwards.
Intuitively, the action-value function quantifies, for a given policy $\pi$, how ``good" the state-action pair $(\bm{s}, \bm{a})$ is in the long run.

Let $Q^*(\cdot, \cdot)$ denote the optimal action-value function---the maximum action-value function over all policies, i.e., $Q^*(\bm{s}, \bm{a}) = \max_{\bm{\pi}} Q^{\bm{\pi}}(\bm{s}, \bm{a})$.
All optimal policies share the same optimal action-value function.
Also, the greedy policy with respect to $Q^*(\bm{s}, \bm{a})$, i.e., $\bm{\pi}^*(\bm{s}) = \argmax_{\bm{a}} Q^*(\bm{s}, \bm{a})$ is an optimal policy.
Then, it follows from \eqref{eq:return} and \eqref{eq:Q_def} that $Q^*(\bm{s}, \bm{a})$ satisfies the following Bellman optimality equation (see, e.g., \cite{sutton2018reinforcement}):
\begin{align} \label{eq:bellman_opt}
	Q^*(\bm{s}, \bm{a}) = \bar{R}(\bm{s}, \bm{a}) + \gamma \sum \limits_{\bm{s'} \in \calS} \calP(\bm{s}' | \bm{s}, \bm{a}) \max_{\bm{a}' \in \calA} Q^*(\bm{s}', \bm{a}').
\end{align}
The MDP is solved if we find $Q^*(\bm{s}, \bm{a})$, and correspondingly, the optimal policy $\bm{\pi}^*$.
It is important to emphasize that \eqref{eq:bellman_opt} is key in solving the MDP.
For ease of notation, in the rest of this paper, we simply write the $Q^*(\bm{s}, \bm{a})$ as $Q(\bm{s}, \bm{a})$.

When both the state and the action sets are finite, the action-value function can be exactly represented in a tabular form that covers all possible pairs $(\bm{s}, \bm{a}) \in \calS \times \calA$.
In this case, if $\calP$ is also known, then the MDP can be solved using, e.g., the so-called policy iteration and value iteration algorithms (see, e.g., \cite{sutton2018reinforcement}).
If $\calP$ is unknown but samples of transitions are available, the MDP can be solved by using RL algorithms such as the Q-learning algorithm (see, e.g., \cite{watkins1992q}).

\subsection{Batch Reinforcement Learning}

When $\calS$ is not finite, conventional Q-learning based approaches require discretization of $\calS$ (see, e.g., \cite{vlachogiannis2004reinforcement} and \cite{xu2012multiagent}).
The discretized state space will better approximate the original state space if a small step size is used in the discretization process, yet the resulting MDP will face the ``curse of dimensionality."
A large step size can alleviate the computational burden caused by the high dimensionality of the state space, but at the cost of potentially degrading performance significantly.

More practically, when the number of elements in $\calS$ is large or $\calS$ is not finite, the action-value function can be approximated by some parametric functions such as linear functions \cite{lagoudakis2003least} and neural networks \cite{mnih2015human}.
Let $\hat{Q}(\cdot, \cdot)$ denote the approximate optimal action-value function.
Using a linear function approximation, $\hat{Q}(\bm{s}, \bm{a})$ can be represented as follows:
\begin{align}
	\hat{Q}(\bm{s}, \bm{a}) = \bm{w}^\top \bm{\phi}(\bm{s}, \bm{a}),
\end{align}
where $\bm{\phi}: \calS \times \calA \rightarrow \bbR^f$ is a feature mapping for $(\bm{s}, \bm{a})$, which is also referred to as the basis function, and $\bm{w} \in \bbR^f$ is the parameter vector.

A class of stable and data-efficient RL algorithms that can solve an MDP with function approximations are the batch RL algorithms---``batch" in the sense that a set of transition samples are utilized each time---such as the LSPI algorithm \cite{lagoudakis2003least}, which is considered to be the most efficient one in this class.
We next explain the fundamental idea behind the LSPI algorithm.
Let $\calD = \{(\bm{s}, \bm{a}, r, \bm{s}'): \bm{s}, \bm{s}' \in \calS, \bm{a} \in \calA\}$ denote a set (batch) of transition samples obtained via observation or simulation.
The LSPI algorithm finds the best $\bm{w}$ that fits the transition samples in $\calD$ in an iterative manner.
One way to explain the intuition behind the LSPI algorithm is as follows (the readers are referred to \cite{lagoudakis2003least} for a more rigorous development).
Define 
\begin{align}
	g(\bm{w}) = \sum_{(\bm{s}, \bm{a}, r, \bm{s}') \in \calD} ( Q(\bm{s}, \bm{a}) - \bm{w}^\top \bm{\phi}(\bm{s}, \bm{a}) )^2.
\end{align}
Let $\bm{w}_i$ denote the value of $\bm{w}$ that is available at the beginning of iteration $i$.
At iteration $i$, the algorithm finds $\bm{w}_{i+1}$ by solving the following problem:
\begin{align} \label{eq:lspi_1}
	\bm{w}_{i+1} = \argmin_{\bm{w}} g(\bm{w}),
\end{align}
which is an unconstrained optimization problem. 
The solution of \eqref{eq:lspi_1} can be computed by setting the gradient of $g(\cdot)$ to zero as follows:
\begin{align} \label{eq:g_grad}
	\frac{\partial g}{\partial \bm{w}} = -2 \sum_{(\bm{s}, \bm{a}, r, \bm{s}') \in \calD} ( Q(\bm{s}, \bm{a}) - \bm{w}^\top \bm{\phi}(\bm{s}, \bm{a})) \bm{\phi}(\bm{s}, \bm{a}) = \zeros_f.
\end{align}

Note that the true value of $Q(\bm{s}, \bm{a})$ is not known and is substituted by the so-called temporal-difference (TD) target, $r + \gamma \bm{w}^\top \bm{\phi}(\bm{s}', \bm{a}')$, where $\bm{a}' = \argmax_{\bm{a} \in \calA} \bm{w}_i^\top \bm{\phi}(\bm{s}', \bm{a})$ is the optimal action in state $\bm{s}'$ determined based on $\bm{w}_i$.
Note that the TD target is a sample of the right-hand-side (RHS) of $\eqref{eq:bellman_opt}$, which serves as an estimate for the RHS of $\eqref{eq:bellman_opt}$.
We emphasize that despite $Q(\bm{s}, \bm{a})$ being substituted by $r + \gamma \bm{w}^\top \bm{\phi}(\bm{s}', \bm{a}')$, the true optimal action-value function is not a function of $\bm{w}$; therefore, the gradient of $g$ with respect to $\bm{w}$ is taken before the $Q(\bm{s}, \bm{a})$ is approximated by the TD target, which does depends on $\bm{w}$.
Then, after replacing $Q(\bm{s}, \bm{a})$ with the TD target, \eqref{eq:g_grad} has the following closed-form solution:
\begin{align} \label{eq:lspi_2}
	\bm{w}_{i+1} =& \left(\sum_{(\bm{s}, \bm{a}, r, \bm{s}') \in \calD} \bm{\phi}(\bm{s}, \bm{a})(\bm{\phi}(\bm{s}, \bm{a}) - \gamma \bm{\phi}(\bm{s}', \bm{a}'))^\top \right)^{-1} \nonumber \\
	& \times \sum_{(\bm{s}, \bm{a}, r, \bm{s}') \in \calD} \bm{\phi}(\bm{s}, \bm{a}) r.
\end{align}

Intuitively, at each iteration, the LSPI algorithm finds the $\bm{w}$ that minimizes the mean squared error between the TD target and $\hat{Q}(\bm{s}, \bm{a})$ over all transition samples in $\calD$.
This process is repeated until change of $\bm{w}$, defined as $\norm{\bm{w}_{i+1} - \bm{w}_i}$, where $\norm{\cdot}$ denotes the $L_2$-norm, becomes smaller than a threshold $\varepsilon$, upon which the algorithm is considered to have converged.

The LSPI algorithm has the following three nice properties.
First, linear functions are used to approximate the optimal action-value function, which allows the algorithm to handle MDPs with high-dimensional or continuous state spaces.
Second, at each iteration, a batch of transition samples is used to update the vector $\bm{w}$ parameterizing  $\hat{Q}(\cdot, \cdot)$, and these samples are reused at each iteration, thus increasing data efficiency.
Third, the optimal parameter vector is found by solving a least-squares problem, resulting in a stable algorithm.
We refer interested readers to \cite{lagoudakis2003least} for more details on the convergence and performance guarantee of the LSPI algorithm.

\section{Optimal Tap Setting Problem as An MDP} \label{sec:formulation}

In this section, we formulate the optimal tap setting problem as an MDP as follows:

\subsubsection{State space}

Define the squared voltage magnitudes at all buses but bus $0$ and the tap ratios as the state, i.e., $\bm{s}~=~(\bm{t}, \bm{v})$, which has both continuous and discrete variables.
Then, the state space is $\calS \subseteq \calT^{L^t} \times \bbR^N$.

\subsubsection{Action space}

The actions are the LTC tap ratio changes, i.e., $\bm{a} = \Delta \bm{t}$, and the action space is the set of all feasible values of LTC tap ratios, i.e., $\calA = \Delta \calT^{L^t}$.
In the optimal tap setting problem, the action is discrete.
The size of the action space increases exponentially with the number of LTCs.


\subsubsection{Reward function}
The objective of voltage regulation is to minimize the voltage deviation as measured by the $L_2$ norm.
As such, when the system transitions from state $\bm{s}=(\bm{t}, \bm{v})$ into state $\bm{s}'=(\bm{t}', \bm{v}')$ after taking action $\bm{a} = \Delta \bm{t} := \bm{t}' - \bm{t}$, the reward is computed by the following function:
\begin{align} \label{eq:reward}
	\calR(\bm{s}, \bm{a}, \bm{s}') = -\frac{1}{N}\norm{\bm{v}' - \bm{v}^{\star}}.
\end{align}

\subsubsection{Transition model}
To derive the transition model $\calP$, note that it follows from \eqref{eq:LDF-Vec-LTC} that
\begin{align} \label{eq:state_trans}
    \bm{v}' =& (\bm{M}(\bm{t}')^\top)^{-1} (\bm{\xi} + \bm{M}(\bm{t})^\top \bm{v} + \bm{m}(\bm{t}) v_0 - \bm{m}(\bm{t}') v_0),
\end{align}
where $\bm{\xi} = 2\diag{\bm{r}} \bm{M}^{-1} (\bm{p}' - \bm{p}) + 2\diag{\bm{x}} \bm{M}^{-1} (\bm{q}' - \bm{q})$, and $\bm{p}'$ and $\bm{q}'$ are active and reactive power injections that results into $\bm{v}'$, respectively.
Then, the transition model $\calP(\bm{s}' | \bm{s}, \bm{a})$ can be derived from the probability density function (pdf) of $(\bm{v}' | \bm{v}, \bm{t}, \Delta \bm{t})$, which can be further computed from the pdf of $(\bm{\xi} | \bm{v}, \bm{t}, \Delta \bm{t})$.
However, under Assumptions \textbf{A1} and \textbf{A2}, the line parameters as well as the probability distributions of active and reactive power injections are unknown; thus, the transition model is not known a priori.
Therefore, we need to resort to RL algorithms that do not require an explicit transition model to solve the MDP.

\section{Optimal Tap Setting Algorithm} \label{sec:algo}

In this section, we propose an optimal tap setting algorithm, which consists of a transition generating algorithm that can generate samples of transitions in $\calD$, and an LSPI-based sequential learning algorithm to solve the MDP.
Implementation details such as the feature selection are also discussed.

\subsection{Overview}

\begin{figure}[!t]
\centering
\includegraphics[width=3.2in]{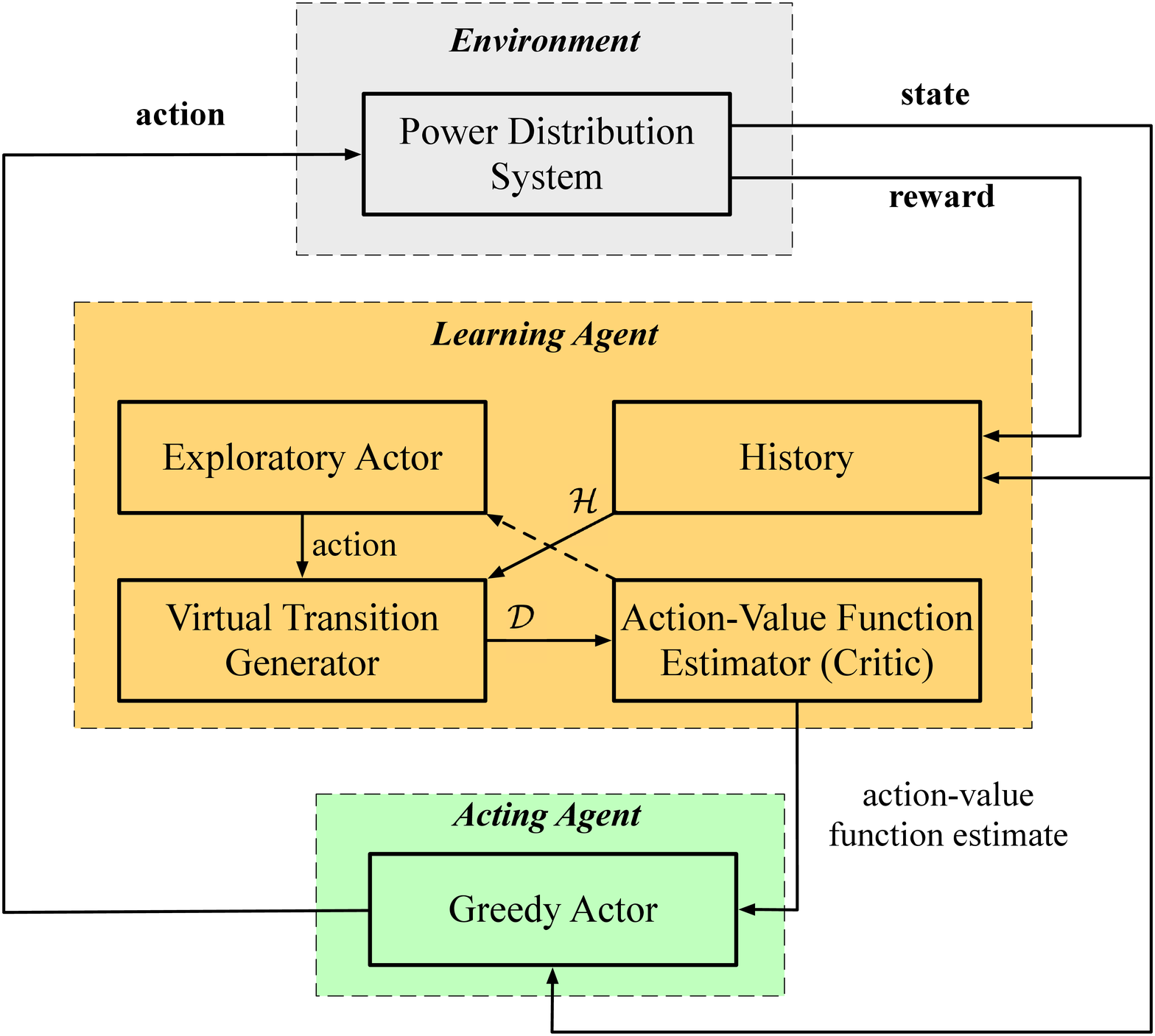}
\caption{The batch RL based framework for optimal tap setting. (Dotted line indicates the critic is optional for the exploratory actor.)}
\label{fig:framework}
\end{figure}

The overall structure of the optimal tap setting framework is illustrated in Fig. \ref{fig:framework}.
The framework consists of an environment that is the power distribution system, a learning agent that learns the action-value function from a set of transition samples, and an acting agent that determines the optimal action from the action-value function.
Define the history to be the sequence of states, actions, and rewards, and denote it by $\calH$, i.e., $\calH = \{\bm{s}_0, \bm{a}_0, r_0, \bm{s}_1, \bm{a}_1, r_1, \cdots\}$.
Specifically, the learning agent will use the elements in the set $\calH$ together with a virtual transition generator to generate a set of transition samples $\calD$ according to some exploratory behavior defined in the exploratory actor.
The set of transition samples in $\calD$ is then used by the action-value function estimator---also referred to as the critic---to fit an approximate action-value function using the LSPI algorithm described earlier.
The learning agent, which has a copy of the up-to-date approximate action-value function from the learning agent, finds a greedy action for the current state and instructs the LTCs to follow it.

Note that the learning of the action-value function can be done offline by the learning agent, which is capable of exploring various system conditions through the virtual transition generator based on the history $\calH$, yet without directly interacting with the power distribution system.
This avoids jeopardizing system operational reliability, which is a major concern when applying RL algorithms to power system applications \cite{glavic2017reinforcement}.

\subsection{Virtual Transition Generator}

The LSPI algorithm (as well as all other RL algorithms) require adequate transition samples that spread over the state and action spaces $\calS \times \calA$.
However, this is challenging in power systems since the system operational reliability might be jeopardized when exploring randomly.
One way to work around this issue is to use simulation models, rather than the physical system, to generate virtual transitions.
To this end, we develop a data-driven virtual transition generator that simulates transitions without any knowledge of the active and reactive power injections (neither measurements nor probability distributions) or the line parameters.

The fundamental idea is the following.
For a transition sample $(\bm{s}, \bm{a}^\dag, r^\dag, \bm{s}^\dag=(\bm{t}^\dag, \bm{v}^\dag))$ that is obtained from $\calH$, the virtual transition generator generates a new transition sample $(\bm{s}, \bm{a}^\ddag, r^\ddag, \bm{s}^\ddag=(\bm{t}^\ddag, \bm{v}^\ddag))$, where $\bm{a}^\ddag$ is determined from $\bm{s}$ according to some exploration policy (to be defined later) that aims to explore the state and action spaces.
Replacing $\bm{a}^\dag$ in the first transition sample with $\bm{a}^\ddag$, the voltage magnitudes will change accordingly.
Assume the same transition of the power injections in these two samples, then the RHS of \eqref{eq:LDF-Vec-LTC} does not change.
Thus, $\bm{v}^\ddag$ can be readily computed from $\bm{v}^\dag$ by solving the following set of linear equations:
\begin{equation} \label{eq:tap-change}
	\bm{M}(\bm{t}^\ddag)^\top \bm{v}^\ddag + \bm{m}(\bm{t}^\ddag) v_0 = \bm{M}(\bm{t}^\dag)^\top \bm{v}^\dag + \bm{m}(\bm{t}^\dag) v_0.
\end{equation}
Since the only unknown in \eqref{eq:tap-change} is $\bm{v}^\ddag \in \bbR$ and $M(\bm{t}^\ddag) \in \bbR^{N \times N}$ is invertible, we can solve for $\bm{v}^\ddag$ as follows:
\begin{align} \label{eq:veg}
	\bm{v}^\ddag = (\bm{M}(\bm{t}^\ddag)^\top)^{-1} (\bm{M}(\bm{t}^\dag)^\top \bm{v}^\dag + \bm{m}(\bm{t}^\dag) v_0 - \bm{m}(\bm{t}^\ddag) v_0).
\end{align}
For ease of notation, we simply write \eqref{eq:veg} as 
\begin{align} \label{eq:est_v}
	\bm{v}^\ddag = \varphi(\bm{v}^\dag, \bm{t}^\dag, \bm{t}^\ddag).
\end{align}
This nice property allows us to estimate the new values of voltage magnitudes when the tap positions change without knowing the exact values of power injections and line parameters.
The virtual transition generating procedure is summarized in Algorithm \ref{algo:transition-generating}.

\begin{algorithm}[!t]
    \SetAlgoLined
    \DontPrintSemicolon
	\KwData{$\calH$, $D$, $\bm{v}^\star$, exploration policy}
    \KwResult{$\calD$}
    Initialize $\calD \leftarrow \varnothing$\;
    \For{$d=1, \cdots, D$}{
    	Choose a transition sample $(\bm{s}, \bm{a}^\dag, r^\dag, \bm{s}^\dag=(\bm{t}^\dag, \bm{v}^\dag))$ from $\calH$\;
        Select $\bm{a}^\ddag$ according to exploration policy and set $\bm{t}^\ddag = \bm{t}^\dag + \bm{a}^\ddag$\;
        Estimate $\bm{v}^\ddag$ following $\bm{a}^\ddag$ as $\bm{v}^\ddag = \varphi(\bm{v}^\dag, \bm{t}^\dag, \bm{t}^\ddag)$\;
		Compute the reward by $r^\ddag = -\frac{1}{N}\norm{\bm{v}^\ddag - \bm{v}^{\star}}$\;
		Add $(\bm{s}, \bm{a}^\ddag, r^\ddag, \bm{s}^\ddag=(\bm{t}^\ddag, \bm{v}^\ddag))$ to $\calD$\;
    }
\label{algo:transition-generating}
\caption{Virtual transition Generating}
\end{algorithm}

\subsection{LSPI-based Sequential Action-Value Function Learning} \label{sec:lspi_seq}

Given the transition sample set $\calD$, we can now develop a learning algorithm for $\hat{Q}(\bm{s}, \bm{a})$ based on the LSPI algorithm.
While the LSPI is very efficient when the action space is relatively small, it becomes computationally intractable when the action space is large, since the number of unknown parameters in the approximate action-value function is typically proportional to $|\calA|$, which increases exponentially with the number of LTCs.
To overcome the ``curse of dimensionality" that results from the size of the action space, we propose an LSPI-based sequential learning algorithm to learn the action-value function.

The key idea is the following. Instead of learning an approximate optimal action-value function for the action vector $\bm{a}$, we learn a separate approximate action-value function for each component of $\bm{a}$. 
To be more specific, for each LTC $l$, $l = 1, \cdots, L^t$, we learn an approximate optimal action-value function $\hat{Q}^{(l)}(\bm{s}, a^{(l)}) = \bm{\phi}^{(l)}(\bm{s}, a^{(l)})^\top \bm{w}^{(l)}$, where $a^{(l)}$ is the $l^{\text{th}}$ component of $\bm{a}$, $\bm{\phi}^{(l)}(\cdot, \cdot)$ is a feature mapping from $\calS \times \Delta \calT$ to $\bbR^f$.
During the learning process of $\bm{w}^{(l)}$, the rest of the LTCs are assumed to behave greedily according to their own approximate optimal action-value function.
To achieve this, we design the following exploration policy to generate the virtual transition samples $\calD$ used when learning $\bm{w}^{(l)}$ for LTC $l$.
In the exploration step in Algorithm \ref{algo:transition-generating}, the tap ratio change of LTC $l$ is selected uniformly in $\Delta \calT$ (uniform exploration), while those of others are selected greedily with respect to the up-to-date $\hat{Q}^{(l)}(\cdot, \cdot)$ (greedy exploration).
Then, the LSPI algorithm detailed in Algorithm \ref{algo:lspi}, where $c$ is a small positive pre-condition number and $\bm{w}_1^{(l)}$ is the initial value for the parameter vector, is applied to learn $\bm{w}^{(l)}$.
This procedure is repeated in a round-robin fashion for all LTCs for $J$ iterations, in each of which $\bm{w}_1^{(l)}$ is set to the up-to-date $\bm{w}^{(l)}$ learned in the previous iteration or chosen if it is in the first iteration.
The value of $J$ is set to $1$ if there is only one LTC and is increased slightly when there are more LTCs.
Note that a new set of transitions $\calD$ is generated when learning $\bm{w}^{(l)}$ for different LTCs at each iteration.
Using this sequential learning algorithm, the total number of unknowns is then proportional to $L^t |\Delta \calT|$, which is far fewer compared to $|\Delta \calT^{L^t}|$ as in the case where the approximate optimal action-value function for the entire action vector, $\bm{a}$, is learned.

A critical step in implementing the LSPI algorithm is constructing features from the state-action pair $(\bm{s}, a^{(l)})$ for LTC $l$; we use radial basis function (RBFs) to this end.
The feature vector for a state-action pair $(\bm{s}, a^{(l)})$, i.e., $\bm{\phi}^{(l)}(\bm{s}, a^{(l)})$, is a vector in $\bbR^f$, where $f = (\kappa+1)\times |\Delta \calT|$ and $\kappa$ is a positive integer.
$\bm{\phi}^{(l)}(\bm{s}, a^{(l)})$ has $|\Delta \calT|$ segments, each one of length $\kappa+1$ corresponding to a tap change in $\Delta \calT$, i.e, $\bm{\phi}^{(l)}(\bm{s}, a^{(l)}) = [\bm{\psi}_1^\top, \cdots, \bm{\psi}_{|\Delta \calT|}^\top]^\top$, where $\bm{\psi}_i \in \bbR^{\kappa + 1}, i = 1, \cdots, |\Delta \calT|$.
Specifically, for $\bm{s} = (\bm{t}, \bm{v})$ and $a^{(l)}$ being the $i^{\text{th}}$ tap change in $\Delta \calT$, $\bm{\psi}_j = \zeros_{\kappa + 1}$ for $j \neq i$, and ${\bm{\psi}_i = [1, e^{-\frac{\norm{\tdbdv - \bar{\bm{v}}_1}}{\sigma^2}}, \cdots, e^{-\frac{\norm{\tdbdv - \bar{\bm{v}}_\kappa}}{\sigma^2}}]^\top}$, where $\sigma > 0$, $\tdbdv = \varphi(\bm{v}, \bm{t}, \tdbdt)$ with $\tdbdt$ being obtained by replacing the $l^{\text{th}}$ entry in $\bm{t}$ with $1$, and $\bar{\bm{v}}_i$, $i=1,\cdots, \kappa$ are pre-specified constant vectors in $\bbR^N$ referred to as the RBF centers.
The action $a^{(l)}$ only determines which segment will be non-zero.
Thus, $\tdbdv$ is indeed the squared voltage magnitudes under the same power injections if the tap of LTC $l$ is at position $0$.
Each RBF computes the distance between $\bm{v}'$ and some pre-specified squared voltage magnitudes.

\begin{algorithm}[!t]
    \SetAlgoLined
    \DontPrintSemicolon
    \KwData{$l$, $\calD$, $\bm{\phi}$, $\gamma$, $\varepsilon$, $c$, $\bm{w}^{(l)}_1$}
    \KwResult{$\bm{w}^{(l)}$}
    Initialize $\bm{w}^{(l)}_0 = \zeros_f$ and $i = 1$\;
    \While{$ \norm{\bm{w}^{(l)}_i - \bm{w}^{(l)}_{i-1}} > \varepsilon$ or $i = 1$}{
    	Initialize $\bm{B}_0 = c \bm{I}_{f\times f}$ and $\bm{b}_0 = \zeros_f$, set $j=1$\;
        \For{$(\bm{s}, \bm{a}, r, \bm{s}') \in \calD$}{
			${a^{(l)'} = \argmax_{a \in \Delta \calT} \bm{\phi}(\bm{s}', a)^\top \bm{w}_i^{(l)}}$\;
        	${\bm{B}_j = \bm{B}_{j-1} + \bm{\phi}(\bm{s}, a^{(l)}) ( \bm{\phi}(\bm{s}, a^{(l)}) - \gamma \bm{\phi}(\bm{s}', a^{(l)'}) )^\top}$\;\vspace{-0.1in}
        	${\bm{b}_j = \bm{b}_{j-1} + \bm{\phi}(\bm{s}, a^{(l)}) r}$\;
        	Increase $j$ by $1$\;
        }
        $\bm{w}^{(l)}_{i+1} = \bm{B}_{|\calD|}^{-1} \bm{b}_{|\calD|}$, increase $i$ by $1$\;
    }
\label{algo:lspi}
\caption{LSPI for Single LTC}
\end{algorithm}

\subsection{Tap Setting Algorithm}

\begin{algorithm}[!t]
    \SetAlgoLined
    \DontPrintSemicolon
	\KwData{$\bm{\phi}$, $K$, $J$, $\epsilon$}
    \For{$k=1, 2, \cdots$}{
    	Obtain $r_{k-1}$ and $\bm{s}_k$, and add them into $\calH$\;
		\If{$k \mod K = 0$}{
			\For{$j=1, \cdots, J$}{
	    		\For{$l = 1, \cdots, L^t$}{
	    			Run Algo. \ref{algo:transition-generating} to generate $\calD$ using uniform exploration for LTC $l$ and greedy exploration for other LTCs\;
		    		Run Algo. \ref{algo:lspi} with $\bm{w}^{(l)}_1$ set to the current $\bm{w}^{(l)}$\;
	    		}\;\vspace{-0.2in}
	    	}
		}
		\For{$l = 1, \cdots, L^t$}{
		Set $a_k^{(l)} = \argmax \limits_{a \in \Delta \calT} \bm{\phi}(\bm{s}_k, a)^\top \bm{w}^{(l)}$ if $\max \limits_{a \in \Delta \calT} \bm{\phi}(\bm{s}_k, a)^\top \bm{w}^{(l)} - \bm{\phi}(\bm{s}_k, a_{k-1}^{(l)})^\top \bm{w}^{(l)} > \epsilon$\;
		Set $a_k^{(l)} = a_{k-1}^{(l)}$ otherwise\;
		}\;\vspace{-0.2in}
    	Add $\bm{a}_k$ to $\calH$ and adjust tap ratios based on $\bm{a}_k$\;
    }
\label{algo:optimal-tap-setting}
\caption{Optimal Tap Setting}
\end{algorithm}

\begin{figure}[!t]
\centering
\includegraphics[width=3.in]{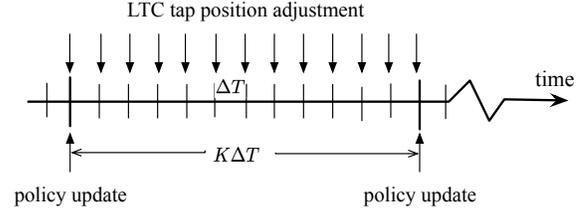}
\caption{Timeline for LTC tap setting.}
\label{fig:timeline}
\end{figure}

The tap setting algorithm, the timeline of which is illustrated in Fig. \ref{fig:timeline}, works as follows.
At time instant $k$, a new state $\bm{s}_k$ as well as the reward following the action $\bm{a}_{k-1}$, $r_{k-1}$, is observed.
Let $\Delta T$ denote the time ellapsed between two time instants.
Every $K$ time instants, i.e., every $K \Delta T$ units of time, $\bm{w}^{(l)}$, is updated by the learning agent by executing the LSPI-based sequential learning algorithm described in Section \ref{sec:lspi_seq}.
The acting agent then finds a greedy action for the current state $\bm{s}_k$ and sends it to the LTCs.
In order to reduce the wear and tear on the LTCs, the greedy action for the current state $\bm{s}_k$ is chosen only if the difference between the action-value resulting from the greedy action, i.e., $\max \limits_{a \in \Delta \calT} \bm{\phi}(\bm{s}_k, a)^\top \bm{w}^{(l)}$, and that resulting from the previous action, i.e., $\bm{\phi}(\bm{s}_k, a_{k-1}^{(l)})^\top \bm{w}^{(l)}$, is larger than a threshold $\epsilon$. 
Otherwise, the tap positions do not change.
The above procedure is summarized in Algorithm~\ref{algo:optimal-tap-setting}.

\section{Numerical Simulation} \label{sec:simu}

In this section, we apply the proposed methodology to the IEEE 13-bus and 123-bus test feeders from \cite{test_feeder}.

\subsection{Simulation Setup}

The power injections for both these two test feeders are constructed based on historical hourly active power load data from a residential building in San Diego over one year \cite{load_data}.
Specifically, the historical hourly active power load data are first scaled up so that the maximum system total active power load over that year for the IEEE 13-bus and 123-bus distribution test feeders are $6.15$~MW and $12.3$~MW, respectively.
These numbers are chosen so that the resulting voltage magnitudes fall outside of the desired range at some time instants.
Then, the time granularity of the scaled system total active power load is increased to $5$ minutes through a linear interpolation.
Each value in the resulting five-minute system total active power load data time series is further multiplied by a normally distributed variable, the mean and standard deviation of which is $1$ and $0.02$, respectively.
The active power load profile at each bus is constructed by pseudo-randomly redistributing the system total active power load among all load buses.
Each load bus is assumed to have a constant power factor of $0.95$.
While only load variation is considered in the simulation, the proposed methodology can be directly applied to the case with renewable-based resources, which can be modeled as negative loads.



We first verify the accuracy of the virtual transition generating algorithm.
Specifically, assume the voltage magnitudes are known for some unknown power injections under a known tap ratio of $1$.
Then, when the tap ratio changes, we compute the true voltage magnitudes under the new tap ratio, denoted by $\bm{v}$, by solving the full ac power flow problem, and the estimated voltage magnitudes under new tap ratio, denoted by $\hat{\bm{v}}$, via \eqref{eq:est_v}.
Simulation results indicate that the maximum absolute difference between the true and the estimated voltage magnitude, i.e., $\norm{\bm{v} - \hat{\bm{v}}}_\infty$, is smaller than $0.001$ p.u., which is accurate enough for the application of voltage regulation addressed in this paper.


\subsection{Case Study on the IEEE 13-bus Test Feeder}

Assume $\bm{v}^{\star} = \ones_N$, where $\ones_N$ is an all-ones vector in $\bbR^N$
In the simulation, $21$ RBF centers are used, i.e., $\kappa = 21$. 
Specifically, $\bar{\bm{v}}_i = (0.895 + 0.005i)^2 \times \ones_N$, $i=1, \cdots, 21$.
The duration between two time instants is $\Delta T = 5$~min.
The policy is updated every $2$ hours, i.e., $K=24$.
In each update, actual transition samples are chosen from the history over the same time interval in the previous $5$ days, which are part of $\calH$, and new actions are chosen according to the exploration policy described in Section \ref{sec:lspi_seq}.
A total number of $D=6000$ virtual transitions are generated using Algorithm \ref{algo:transition-generating}.
Since this test feeder only has one LTC, there is no need to sequentially update the approximate action-value function, so we set $J=1$.
Other parameters are chosen as follows: $\gamma=0.9$, $\varepsilon=1\times 10^{-5}$, $\epsilon=1\times 10^{-4}$, $c=0.1$, and $\sigma=1$.

Assuming complete and perfect knowledge on the system parameters as well as active and reactive power injections for all time instants, we can find the optimal tap position that results in the highest reward by exhaustively searching the action space, i.e., all feasible tap ratios, at each time instant.
It is important to point out that, in practice, the exhaustive search approach is infeasible since we do not have the necessary information, and not practical due to the high computational burden.
Results obtained by the exhaustive search approach and the conventional tap setting scheme (see, e.g., \cite{kundur1994power}), in which the taps are adjusted only when the voltage magnitudes exceed a desired range, e.g., $[0.9, 1.1]$~p.u., are used to benchmark the proposed methodology.

%

\begin{figure}[!t]
\centering
\includegraphics[width=3.5in]{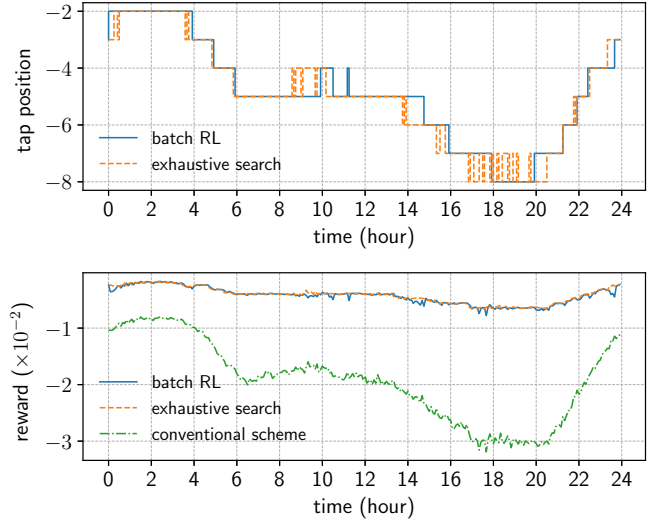}
\caption{Tap positions and rewards for IEEE 13-bus test feeder.}
\label{fig:tap_reward_13_bus}
\end{figure}

\begin{figure}[!t]
\centering
\includegraphics[width=3.5in]{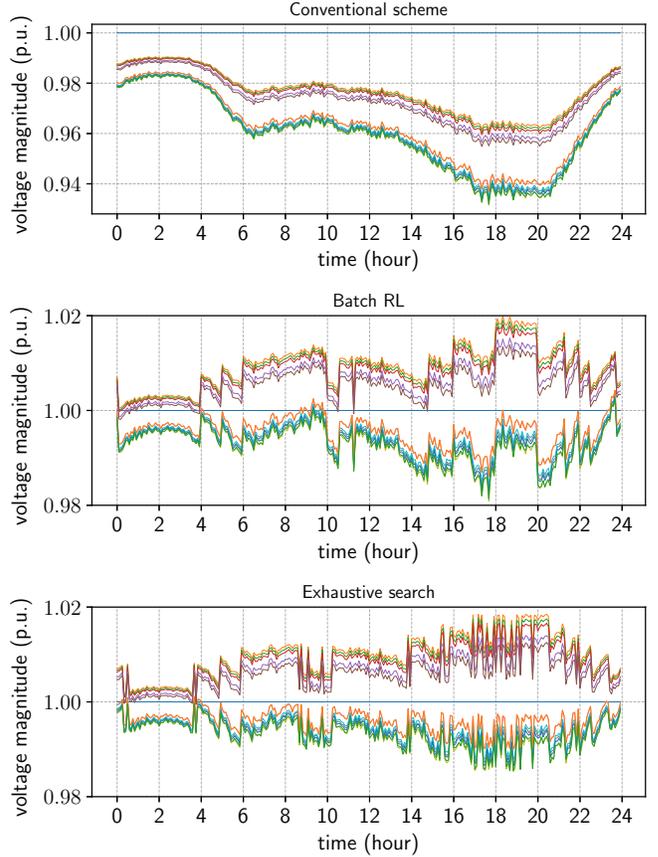}
\caption{Voltage magnitude profiles of IEEE 13-bus test feeder.}
\label{fig:voltage_13_bus}
\end{figure}

Figure \ref{fig:tap_reward_13_bus} shows the tap positions (top panel) and the rewards (bottom panel) under different approaches.
The rewards resulted from these two approaches are very close.
The daily mean reward, i.e., $\rho = \frac{1}{288} \sum_{k=1}^{288} r_k$, where $r_k$ is the reward at time instant $k$ as defined in \eqref{eq:reward}, obtained by the batch RL approach and the exhaustive search approach is $\rho = -4.279\times 10^{-3}$ and $\rho = -4.156\times 10^{-3}$, respectively, while that under the conventional scheme is $\rho = -19.26 \times 10^{-3}$.
The tap positions under the batch RL approach and the exhaustive search approach are aligned during most of the time during the day.
Note that the tap position under the conventional scheme remains at $0$ since the voltage magnitudes are within $[0.9, 1.1]$~p.u., and is not plotted.
Figure \ref{fig:voltage_13_bus} shows the voltage magnitude profiles under the different tap setting algorithms.
The voltage magnitude profiles under the proposed batch RL approach (see Fig. \ref{fig:voltage_13_bus}, center panel) are quite similar to those obtained via the exhaustive search approach (see Fig. \ref{fig:voltage_13_bus}, bottom panel), both result in a higher daily mean reward than that resulted from the conventional scheme (see Fig. \ref{fig:voltage_13_bus}, top panel).
We also would like to point out that Algorithm \ref{algo:lspi} typically converges within $5$ iterations in less than $20$ seconds, and the batch RL approach is faster than the exhaustive search approach by several orders of magnitude.

%
%


\subsection{Case Study on the IEEE 123-bus Test Feeder}

We next test the proposed methodology on the IEEE 123-bus test feeder.
In the results for the IEEE 13-bus test feeder reported earlier, while the LTC has $33$ tap positions, only a small portion of them is actually used.
This motivates us to further reduce the action space by narrowing the action space to a smaller range.
Specifically, we can estimate the voltage magnitudes under various power injections and LTC tap positions using \eqref{eq:est_v}.
After ruling out tap positions under which the voltage magnitudes will exceed the desired range, we eventually allow $9$ positions, from $-8$ to $0$, for two LTCs, and $5$ positions, from $0$ to $5$, for the other two LTCs.
Here, $\kappa = 11$ RBF centers are used.
Specifically, $\bar{\bm{v}}_i = (0.94 + 0.01i)^2 \times \ones_N$ for all LTCs except for the one near the substation, for which $\bar{\bm{v}}_i = (0.89 + 0.01i)^2 \times \ones_N$, $i=1, 2, \cdots, 11$.
A total number of $D=3600$ virtual transitions are generated in a similar manner as in the IEEE 13-bus test feeder case.
The number of iterations in the LSPI-based sequential learning algorithm is set to $J=3$.
Other parameters are the same as in the IEEE 13-bus test feeder case.


\begin{figure}[!t]
\centering
\includegraphics[width=3.5in]{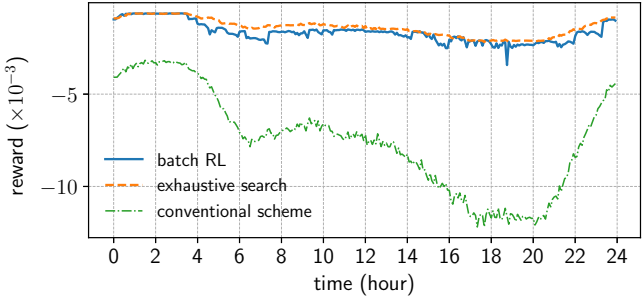}
\caption{Rewards for IEEE 123-bus test feeder.}
\label{fig:reward_123_bus}
\end{figure}


Figure \ref{fig:reward_123_bus} shows the rewards under the batch RL approach and the exhaustive search.
The daily mean reward obtained by the batch RL approach and the exhaustive search approach is $\rho = -1.646\times 10^{-3}$ and $\rho = -1.402\times 10^{-3}$, respectively, while that under the conventional scheme is $\rho = -7.513 \times 10^{-3}$.
Due to the space limitation, other simulation results such as voltage profiles are not presented.


\color{black}

\section{Concluding Remarks} \label{sec:con}

In this paper, we formulate the optimal tap setting problem of LTCs in power distribution systems as an MDP and propose a batch RL algorithm to solve it.
To obtain adequate state-action samples, we develop a virtual transition generator that estimates the voltage magnitudes under different tap settings.
To circumvent the ``curse of dimensionality", we proposed an LSPI-based sequential learning algorithm to learn an action-value function for each LTC, based on which the optimal tap positions can be determined directly.
The proposed algorithm can find the policy that determines the optimal tap positions that minimize the voltage deviation across the system, based only on voltage magnitude measurements and network topology information, which makes it more desirable for implementation in practice.
Numerical simulation on the IEEE 13- and 123-bus test feeders validated the effectiveness of the proposed methodology.

\bibliographystyle{IEEEtran}
\bibliography{LTC}

\end{document}